\documentclass[11pt]{article}
\usepackage{amssymb,amsfonts,amsmath,latexsym,epsf,tikz,url}

\newtheorem{theorem}{Theorem}[section]
\newtheorem{proposition}[theorem]{Proposition}

\newtheorem{corollary}[theorem]{Corollary}

\newcommand{\proof}{\noindent{\bf Proof.\ }}
\newcommand{\qed}{\hfill $\square$\medskip}

\begin{document}

\title{On the structure of dominating graphs}

\author{
Saeid Alikhani $^{a,}$\footnote{Corresponding author}
\and
Davood Fatehi $^{a}$
\and
Sandi Klav\v zar $^{b,c,d}$
}

\date{\today}

\maketitle

\begin{center}
$^a$ Department of Mathematics, Yazd University, 89195-741, Yazd, Iran\\
{\tt alikhani@yazd.ac.ir, davidfatehi@yahoo.com}\\

\medskip
$^b$ Faculty of Mathematics and Physics, University of Ljubljana, Slovenia\\
{\tt sandi.klavzar@fmf.uni-lj.si}\\
\medskip

$^c$ Faculty of Natural Sciences and Mathematics, University of Maribor, Slovenia\\
\medskip

$^d$ Institute of Mathematics, Physics and Mechanics, Ljubljana, Slovenia\\
\end{center}


\begin{abstract}
The $k$-dominating graph $D_k(G)$ of a graph $G$ is defined on the vertex set consisting of dominating sets of $G$ with cardinality at most $k$, two such sets being adjacent if they differ by either adding or deleting a single vertex. A graph is a  dominating graph if it is isomorphic to $D_k(G)$ for some graph $G$ and some positive integer $k$.  Answering a question of Haas and Seyffarth for graphs without isolates, it is proved that if $G$ is  such a graph of order $n\ge 2$ and with $G\cong D_k(G)$, then $k=2$ and $G=K_{1,n-1}$ for some $n\ge 4$. It is also proved that for a given $r$ there exist only a finite number of $r$-regular, connected dominating graphs of connected graphs.  In particular, $C_6$ and $C_8$ are the only dominating graphs in the class of cycles. Some results on the order of dominating graphs  are also obtained.
\end{abstract}


\medskip
\noindent{\bf AMS Subj.\ Class.:} 05C69, 05C62

\section{Introduction}

Let $S$ and $S'$ be dominating sets of a graph $G$ of order at most $k$, where $k$ is a given threshold. Then the {\em dominating set reconfiguration (DSR) problem} asks whether there exists a sequence of dominating sets of $G$ starting with $S$ and ending with $S'$, such that each dominating set in the sequence is of order at most $k$ and can be obtained from the previous one by either adding or deleting exactly one vertex. The problem is PSPACE-complete even for planar graphs, bounded bandwidth graphs, split graphs, and bipartite graphs, while on the positive side it can be solved in linear time for cographs, trees, and interval graphs~\cite{haddadan-2015}.

 The DSR problem naturally leads to the concept of the $k$-dominating graph introduced by Haas and Seyffarth~\cite{haas-2014} as follows. If $G$ is a graph  and $k$ a positive integer, then the {\em $k$-dominating graph} $D_k(G)$ of $G$ is the graph whose vertices correspond to the dominating sets of $G$ that have cardinality at most $k$, two vertices of $D_k(G)$ being adjacent if and only if the corresponding dominating sets of $G$ differ by either adding or deleting a single vertex. (A similar concept is the one of $\gamma$-graphs in which only minimum dominating sets are considered as vertices of the derived graph~\cite{fricke-2011}.)  Now, the DSR problem simply asks whether given two vertices of $D_k(G)$ belong to the same connected component of $D_k(G)$.  Besides with the DSR problem, the $k$-dominating graphs were further motivated by similar studies of graph colorings and by a general goal to further understand the relationship between dominating sets of a graph. 

 It follows from the above discussion that a fundamental problem about $k$-dominating graphs is to determine conditions which ensure that $D_k(G)$ is connected. This problem was the central theme of the seminal paper~\cite{haas-2014}. It is interesting to observe that the connectedness of $D_k(G)$ does not guarantee the connectedness of $D_{k+1}(G)$. For instance, $D_k(K_{1,n-1})$ ($n\ge 4$) is connected for any $1\le k\le n-2$, but $D_{n-1}(K_{1,n-1})$ is not connected. For the latter fact note that $\Gamma(K_{1,n-1}) = n-1$ and that in general $D_{\Gamma(G)}(G)$ is not connected. (Here $\Gamma(G)$ is the upper domination number of $G$, that is, the maximum cardinality of a minimal dominating set of $G$.) On the other hand, Haas and Seyffarth proved that if $G$ has at least two disjoint edges and $k\geq \min\{n-1, \Gamma(G) + \gamma(G)\}$, then $D_k(G)$ is connected. Moreover, if $G$ is bipartite or chordal, then $D_{\Gamma(G)+1}(G)$ is always connected.  
The connectivity of dominating graphs was further investigated in~\cite{suzuki-2014} where it was in particular demonstrated that there exists an infinite family of graphs such that $D_{\gamma(G)+1}(G)$ has exponential diameter and that $D_{n - \mu}(G)$ is connected for any graph $G$ of order $n$ and with a matching of size at least $\mu+1$.

 In this paper we continue the study of $k$-dominating graphs and proceed as follows. In the next section, we introduce additional concepts needed, recall some basic properties of $k$-dominating graphs, and add additional results to this list. In Section~\ref{sec:self}, we attack the question from~\cite{haas-2014} where it was observed that $D_2(K_{1,n}) \cong K_{1,n}$ and asked whether there are other graphs $G$ for which $D_k(G) \cong G$ holds. We prove that if $G$ is of order $n\ge 2$ and with $\delta \ge 1$, and if $G\cong D_k(G)$, where $\gamma(G)\le k\le n$, then actually $G\cong K_{1,n-1}$ holds for some $n\ge 4$.  Then, in Section~\ref{sec:regular-finite}, we prove that for any $r\ge 1$ there exists only a finite number of $r$-regular, connected dominating graphs of connected graphs. For $r=2$ we strengthen the result by showing that $C_6$ and $C_8$ are the only such graphs. We also show that among the paths, $P_1$ and $P_3$ are the only dominating graphs of connected graphs. In the final section we present some results on the order of $k$-dominating graphs, while along the way several problems for further study are stated.

\section{Preliminaries}

We use the notation $[n] = \{1,\ldots, n\}$. As usual, $\delta(G)$ and $\Delta(G)$ denote the minimum and the maximum degree of a graph $G$, respectively. The order of a graph $G=(V,E)$ is denoted with $|G|$, that is, $|G| = |V|$, and the disjoint union of graphs $G$ and $H$ is denoted with $G\cup H$. The {\em join} $G+H$ of graphs $G$ and $H$ is obtained from the disjoint union of $G$ and $H$ by connecting any vertex of $G$ with any vertex of $H$. We write $G\cong H$ to say that $G$ and $H$ are isomorphic graphs.

If $G = (V,E)$ is a graph, then $S\subseteq V$ is a {\em dominating set} of $G$ if every vertex in $V - S$ is adjacent to at least one vertex in $S$. The {\it domination number} $\gamma(G)$ of $G$ is the minimum cardinality of a dominating set in $G$. A dominating set of the minimum cardinality is called a {\em $\gamma$-set}. A vertex of $G$ of degree $|G|-1$ is called a {\em dominating vertex} of $G$. For additional concepts from the domination theory see~\cite{haynes-1998}.

We say that a graph is a {\em dominating graph} if it is isomorphic to $D_k(G)$ for some graph $G$ and some positive integer $k$. For example, $C_6$ is a dominating graph because $D_2(K_3)\cong C_6$. In the next result we collect several basic properties about dominating graphs.

\begin{proposition}
\label{prp:basic}
If $G$ is a graph, then the following hold.
\begin{itemize}
\item[(i)] If $\gamma(G)\le k\le |G|$, then $D_k(G)$ is bipartite.
\item[(ii)] $|D_{|G|}(G)|$ is odd and $|D_{|G|-1}(G)|$ is even.
\item[(iii)] If $m$ is odd, $0 < m < 2n$, then there exists a graph $X$ of order $n$ such that $|D_{n}(X)| = m$.
\item[(iv)] If $G$ is connected, then $\Delta(D_{|G|}(G)) = |G|$.
\end{itemize}
\end{proposition}

\proof
(i) Note that $D_n(K_n)$ is isomorphic to the graph obtained from the $n$-cube $Q_n$ by deleting one of its vertices. Since $D_k(G)$ is a subgraph of $D_k(K_n)$ and the latter graph is a subgraph of the bipartite graph $D_n(K_n)$, it follows that $D_k(G)$ is bipartite.

(ii) That the order of $D_{|G|}(G)$ is odd follows immediately from a result of Brouwer, Csorba, and Schrijver~\cite[Theorem 1.1]{brouwer-2009} asserting that the number of dominating sets of a finite graph is odd. As the only dominating set of order $n$ of $G$ is its vertex set, $D_{|G|-1}(G)$ is then of even order.

(iii) This assertion follows from~\cite[Proposition 1.2]{brouwer-2009} which asserts that if $m$ is odd, where $0 < m < 2n$, then there exists a graph of order $n$ that contains precisely $m$ dominating subsets (see also~\cite{alikhani-2013a}).

(iv) As $D_{|G|}(G)$ is a subgraph of $Q_{|G|}$ we infer that $\Delta(D_{|G|}(G)) \le |G|$. On the other hand, since $G$ is connected, any $(|G|-1)$-subset of vertices is a dominating set and adjacent to the whole vertex set in $D_{|G|}(G)$. So $V(G)$ is  of degree $|G|$ in $D_{|G|}(G)$.
\qed

$D_{|G|}(G)$ is not regular unless $G$ is an edge-less graph in which case $D_{|G|}(G)\cong K_1$. Note also that from Proposition~\ref{prp:basic}(i) and (ii) it follows that $D_{|G|}(G)$ is not hamiltonian. On the other hand, the question which $k$-dominating graphs $D_k(G)$ with $k < |G|$ are hamiltonian remains as open problem.

\section{Graphs isomorphic to their dominating graphs}
\label{sec:self}

Haas and Seyffarth~\cite{haas-2014} observed that $D_2(K_{1,n}) \cong K_{1,n}$ and posed the question whether there are other graphs $G$ for which $D_k(G) \cong G$. In the next result we prove that the answer is negative as soon as $G$ has no isolated vertices.

\begin{theorem}
\label{thm:stars}
Let $G$ be a graph of order $n\ge 2$ and with $\delta \ge 1$. If $G\cong D_k(G)$, where $\gamma(G)\le k\le n$, then $k=2$ and $G=K_{1,n-1}$ for some $n\ge 4$.
\end{theorem}

\proof
Assume first that $k=\gamma(G)$. Then $V(D_k(G))$ consists of $\gamma$-sets of $G$ and  hence $D_k(G)$ is an edge-less graph. If $G\cong D_k(G)$, this is only possible when $G=K_1$. As we have assumed that $G$ is of order at least $2$, we may suppose in the rest of the proof that $k\ge \gamma(G)+1$.

Let $V(G) = \{v_1,\ldots, v_n\}$ and set $\gamma = \gamma(G)$. Let $X$ be a $\gamma$-set of $G$, where we may without loss of generality assume that $X = \{v_1,\ldots, v_\gamma\}$.
Assume that $G\cong D_k(G)$ and recall that $k\ge \gamma+1$. Then $X_i = X\cup \{v_i\}$, $\gamma+1\le i\le n$, are dominating sets of $G$ and hence vertices of $D_k(G)$. As they are all of cardinality $\gamma+1$, the vertices $X,X_{\gamma+1}, \ldots, X_n$ induce a $K_{1,n-\gamma}$ in $D_k(G)$ (hence $G$ also contains an induced $K_{1,n-\gamma}$). Let ${\cal Y} = \{Y_1,\ldots, Y_{\gamma-1}\}$ be the remaining vertices of $D_k(G)$. Observe that the vertex $Y_i$, $i\in [\gamma-1]$, is not adjacent to $X$, for otherwise $\{X\} \cup \{Y_j:j\ne i\}$ would be a dominating set of $D_k(G)$, but then (since $G\cong D_k(G)$) we would have a dominating set of $G$ smaller than $\gamma$. Moreover, a vertex $X_i$, $\gamma+1\le i\le n$, can be adjacent to at most one vertex from ${\cal Y}$. Indeed, suppose that, without loss of generality, $X_{\gamma+1}$ is adjacent to $Y_1$ and $Y_2$. Then $\{X, X_{\gamma+1}, Y_{3}, \ldots, Y_{\gamma-1}\}$ is a dominating set  of $D_k(G)$ yielding the same contradiction as above.

 If for some $i\ne j$, $Y_i$ would be adjacent to $Y_j$, then $\{X,Y_1,\ldots Y_{\gamma-1}\}\setminus \{Y_j\}$ would be a dominating set of $D_k(G)$ of size $\gamma - 1$. Hence, since by the theorem's assumption $G$ has no isolated vertices, each $Y_i$ has a neighbor in ${\cal X} = \{X_{\gamma+1}, \ldots, X_n\}$. Since furthermore no $X_j$ is adjacent to two vertices from ${\cal Y}$, we find out that there exists a matching from ${\cal Y}$ to ${\cal X}$. Let $\{X_{i_1},\ldots, X_{{i_\gamma-1}}\}$ be the endpoints of the matching edges which lie in ${\cal X}$. Then $X,X_{i_1},\ldots, X_{{i_\gamma-1}}$ is a dominating set of $G$ of cardinality $\gamma$ and we have the following two $\gamma$-sets of $D_k(G)$:
\begin{itemize}
\item $X,Y_{1},\ldots, Y_{\gamma-1}$\quad and
\item $X,X_{i_1},\ldots, X_{{i_\gamma-1}}$.
\end{itemize}
Suppose that $\gamma \ge 2$. Adding to any of the above two $\gamma$-sets an additional vertex, we get a dominating set of cardinality $\gamma +1$. Since $\gamma \ge 2$, in this way we can construct $2(n-\gamma)-1$ different dominating sets of $G$ of this cardinality. Consequently, $D_k(G)$ contains at least $2 + 2(n-\gamma)-1$ vertices. Since for any graph without isolated vertices $\gamma \le n/2$ holds, it follows that $D_k(G)$ contains at least $n+1$ vertices, a contradiction.

The only case left to consider is $\gamma = 1$. Assume without loss of generality that $v_1$ is a dominating vertex. Suppose that $G$ contains another dominating vertex, say $v_2$, that is, ${\rm deg}(v_1) = {\rm deg}(v_2)=n-1$. Then $\{v_1\}$, $\{v_2\}$, $\{v_1, v_2\}$, $\{v_1,v_i\}$ ($i\ge 3$), and $\{v_2,v_i\}$ ($i\ge 3$), are dominating sets of $G$, hence $|D_2(G)|\ge 2n-1 > n$. Therefore $v_1$ is the unique dominating vertex of $G$. Now, since $D_2(G)$ is of order $n$, its dominating sets of order at most $2$ are $\{v_1\}$ and $\{v_1,v_i\}$ ($i\ge 2$). But then $D_2(G)\cong K_{1,n-1}$ where $n\ge 4$.
\qed

Let $G$ be an arbitrary graph with $\gamma(G)\ge 3$ and consider the join $G+K_1$,  where $V(K_1) = \{x\}$. Clearly, $\gamma(G+K_1) = 1$. Moreover, if  $D$ is a dominating set of $G+K_1$ and $|D|=2$, then (since $\gamma(G)\ge 3$) we must have $x\in D$. It follows that $D_2(G+K_1) \cong K_{1,|G|}$. This example shows that the stars $K_{1,n}$ can be represented as  dominating graphs in many different ways.

\section{Realizability of graphs as dominating graphs}
\label{sec:regular-finite}

Another problem from~\cite{haas-2014} is which graphs are dominating graphs.  The main result of this section asserts that not many regular graphs are such. To state the result, a short preparation is needed.

For any $r\ge 1$ let $c_r$ be a given, fixed constant such that $\gamma(G) \le c_r|G|$ holds for any connected graph $G$ with $\delta(G)=r$. As already observed by Ore~\cite{ore-1962}, if $\delta(G)\ge 1$ then $\gamma(G)\le n/2$,  so that we can set $c_1 = 1/2$. The constant $c_2 =2/5$ was independently obtained in~\cite{blank-1973,mccuaig-1989} (actually, there are seven small graphs: $C_4$, and six graphs on seven vertices, for which the  2/5 bound does not hold); the result $c_3 = 3/8$ is due to Reed~\cite{reed-1996}; $c_4 = 4/11$ is from~\cite{sohn-2009}. For $k\ge 5$ the best known constants $c_k$ were recently developed in~\cite{bujtas-2015+}. To obtain these constants a modification of a method from~\cite{bujtas-2015} was applied, which was in turn developed for the investigation of the domination game~\cite{bresar-2010}.

 Let now $r\ge 1$. Then setting 
$${\cal D}_r = \{H:\ H\ {\rm is\ an}\ r\mbox{-}{\rm regular,\, connected \, dominating \ graph\ of\ a\ connected\ graph}\}\,$$
our result reads as follows. 
 
\begin{theorem}
\label{thm:finite}
Let $r\ge 1$. If $G$ is a connected graph such that for some $k$, $D_k(G)\in {\cal D}_r$, then $|G|\le 2r$. Consequently, $|{\cal D}_r| < \infty$. 
\end{theorem}

\proof
 Let $G$ be a connected graph and suppose that $D_k(G)\cong H$, where $H$ is an $r$-regular,  connected graph and $k$ is a positive integer.  Clearly, $k\ge \gamma(G)+1$, for if $k=\gamma(G)$, then $D_k(G)$ is edgeless.  Let $X$ be a $\gamma$-set of $G$. Then for any vertex $y\notin X$, the set $X_y = X\cup \{y\}$ is a dominating set of order $\gamma(G)+1$. Hence $X_y$ is a neighbor of $X$ in $D_k(G)$. Because there are $|G|-\gamma(G)$ such vertices $y$, we infer that $\deg_{D_k(G)}(X)=|G|-\gamma(G)$. As $D_k(G)\cong H$ and $H$ is $r$-regular, it follows that $r=|G|-\gamma(G)$. Since $\gamma(G) \le c_{\delta(G)}|G|$ we have $|G|-r = \gamma(G) \le c_{\delta(G)}|G|$. By the above mentioned Ore's result, $c_{\delta(G)} \le 1/2$ holds, hence we find out that $|G|-r \le |G|/2$ and thus $|G|\le 2r$. 

By the above it follows that for a given $r$, a graph $H\in {\cal D}_r$ can be realized as a dominating graph only with a graph $G$ of order at most $2r$ (and for some fixed $k\le 2r$). As there are only a finite number of such graphs, $|{\cal D}_r| < \infty$.
\qed

Theorem~\ref{thm:finite} strengthens in the case $r=2$ as follows.  
  
\begin{corollary}\label{coro:cycles}
${\cal D}_2 = \{C_6, C_8\}$. 
\end{corollary}

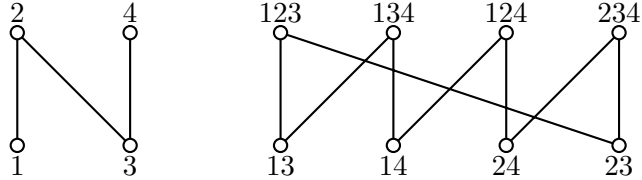
\begin{figure}[ht!]
\begin{center}
\begin{tikzpicture}[scale=0.5,style=thick]
\def\vr{5pt} 

\path (-7,0) coordinate (1);
\path (-7,3) coordinate (2);
\path (-4,0) coordinate (3);
\path (-4,3) coordinate (4);
\path (0,0) coordinate (13);
\path (3,0) coordinate (14);
\path (6,0) coordinate (24);
\path (9,0) coordinate (23);
\path (0,3) coordinate (123);
\path (3,3) coordinate (134);
\path (6,3) coordinate (124);
\path (9,3) coordinate (234);
\draw (1) -- (2) -- (3) -- (4);
\draw (13) -- (134) -- (14) -- (124) -- (24) -- (234) -- (23) -- (123) -- (13);
\draw (1)  [fill=white] circle (\vr);
\draw (2)  [fill=white] circle (\vr);
\draw (3)  [fill=white] circle (\vr);
\draw (4)  [fill=white] circle (\vr);
\draw (13)  [fill=white] circle (\vr);
\draw (14)  [fill=white] circle (\vr);
\draw (23)  [fill=white] circle (\vr);
\draw (24)  [fill=white] circle (\vr);
\draw (123)  [fill=white] circle (\vr);
\draw (134)  [fill=white] circle (\vr);
\draw (124)  [fill=white] circle (\vr);
\draw (234)  [fill=white] circle (\vr);
\draw[below] (1) node {$1$};
\draw[below] (3) node {$3$};
\draw[above] (2) node {$2$};
\draw[above] (4) node {$4$};
\draw[below] (13) node {$13$};
\draw[below] (14) node {$14$};
\draw[below] (23) node {$23$};
\draw[below] (24) node {$24$};
\draw[above] (123) node {$123$};
\draw[above] (134) node {$134$};
\draw[above] (124) node {$124$};
\draw[above] (234) node {$234$};
\end{tikzpicture}
\end{center}
\caption{$P_4$ and $D_3(P_4)$}
\label{fig:C8}
\end{figure}

\proof
 By Theorem~\ref{thm:finite}, $|G|\le 4$ holds if $G$ is a connected graph with $D_k(G)$ isomorphic to a 2-regular connected graph, that is, to a cycle. By inspection we find out that  among connected graphs $G$ of order at most $4$ and among appropriate values $k$, the only favourable cases are $D_2(K_3) \cong C_6$ and $D_3(P_4) \cong C_8$. The fact that $D_3(P_4) \cong C_8$ can be verified using Fig.~\ref{fig:C8}.
\qed

A result parallel to Corollary~\ref{coro:cycles} for paths  reads as follows.

\begin{proposition}\label{prp:paths}
Among the paths, $P_1$ and $P_3$ are the only dominating graphs of connected graphs.
\end{proposition}

\proof By inspection on connected graphs of order at most $4$ the only dominating graphs that are paths are $P_1 \cong D_1(K_1)$ and $P_3 \cong D_2(P_2)$.

Suppose now that $D_k(G)\cong P_m$ holds for some connected graph $G$ with $|G|>4$ and for some $k$ and $m$. Let $X$ be  a $\gamma$-set of $G$. Then either $\deg_{D_k(G)}(X)=1$ or $\deg_{D_k(G)}(X)=2$. Since clearly $k > \gamma(G)$, it follows that either $|G|-\gamma(G) = 1$ or $|G|-\gamma(G) = 2$. But this is not possible since  $|G|>4$.
\qed

In Corollary~\ref{coro:cycles} and in Proposition~\ref{prp:paths} we have considered the dominating graphs that are derived from connected graphs. The following examples indicate that it would be interesting to extend the investigation to disconnected graphs: $D_3(K_2\cup K_2) \cong C_8$ and $D_3(K_2\cup K_1) \cong D_4(K_2\cup K_1\cup K_1) \cong P_3$. Similarly, in Theorem~\ref{thm:finite} we have assumed that the graph $G$ considered has no isolated vertices, hence an extension to graphs that contain isolates could also be interesting.

\section{On the order of dominating graphs}

The {\em domination polynomial} $D(G,x)$ of $G$ is defined as
$$D(G,x)=\sum_{i\ge 0} d(G,i) x^{i}\,,$$
where $d(G,i)$ is the number of dominating sets of $G$ of cardinality $i$. This graph polynomial was introduced in the paper~\cite{alikhani-2014} that appeared in 2014 but numerous other papers on the polynomial appeared earlier. For some very recent developments on the polynomial see~\cite{anthony-2015}. From our perspective, a key information encoded into the domination polynomial is that
$$|D_{|G|}(G)| = D(G,1)\,.$$
For instance, using a result from~\cite{alikhani-2010} asserting that $D(C_1, x) = x$, $D(C_2, x) = x^2 + 2x$, $D(C_3, x) = x^3 + 3x^2 + 3x$, and $D(C_n, x) = x\left( D(C_{n-1}, x) + D(C_{n-2}, x) + D(C_{n-3}, x)\right)$ for $n\ge 4$, we get the following result.

\begin{proposition}
$|D_1(C_1)|=1$, $|D_2(C_2)|=3$, $|D_3(C_3)|=7$, and 
$$|D_n(C_n)|=|D_{n-1}(C_{n-1})|+|D_{n-2}(C_{n-2})|+|D_{n-3}(C_{n-3})|,\ n\geq 4\,.$$
\end{proposition}

We conclude the paper by determining the order of the dominating graph of the join and the corona of two graphs in terms of the invariants of their factors. The join has already been defined, while the {\em corona} $G\circ H$ of graphs $G$ and $H$ is the graph obtained from the disjoint union of $G$ and $|G|$  copies of $H$ by joining the $i$th vertex of $G$ ($1\le i\le |G|$) to every vertex in the $i$-th copy of $H$.

\begin{proposition}
If $G$ and $H$ are graphs, then
\begin{itemize}
\item[(i)] $|D_{|G+H|}(G+H)|=(2^{|G|}-1)(2^{|H|}-1)+|D_{|G|}(G)|+|D_{|H|}(H)|$,
\item[(ii)] $|D_{|G\circ H|}(G\circ H)|=(2^{|H|} +|D_{|H|}(H)|)^{|G|}$.
\end{itemize}
\end{proposition}

\proof
(i) Note first that if $\emptyset \ne D_G \subseteq V(G)$ and $\emptyset \ne D_H \subseteq V(G)$, then $D_G\cup D_H$ is a dominating set of $G+H$. This gives $(2^{|G|}-1)(2^{|H|}-1)$ dominating sets of $G+H$. Assume now that $D$ is a dominating set of $G+H$ with $D\cap V(G)=\emptyset$. Then $D\cap V(H)$ must be a dominating set of $H$, whence there are $|D_{|H|}(H)|$ such dominating sets of $G+H$. Analogously, if $D\cap V(H)=\emptyset$ we get $|D_{|G|}(G)|$ dominating sets of $G+H$.

(ii) Let $D$ be a dominating set of $G\circ H$ and assume that a vertex $x\in V(G)$ does not belong to $D$. If $H_x$ is the copy of $H$ corresponding to $x$, then $D\cap V(H_x)$ is a dominating set of $D_x$. Therefore, if a vertex $x\in V(G)$ is not in $D$, it is dominated by a vertex from $H_x$. It follows that $D$ is a dominating set of $G\circ H$ if and only if $D$ is a dominating set of the graph $(G\circ H) - E(G)$. The latter graph is isomorphic to the disjoint union of $|G|$ copies of the graph $K_1 + H$. By (i), $D_{|H|+1}(K_1+H) = 2^{|H|} +|D_{|H|}(H)|$ and consequently $|D_{|G\circ H|}(G\circ H)|=(2^{|H|} +|D_{|H|}(H)|)^{|G|}$.
\qed


\begin{thebibliography}{99}
\bibitem{alikhani-2013a}
  S.~Alikhani,
  The domination polynomial of a graph at $-1$,
  Graphs Combin.\ 29 (2013) 1175--1181.

\bibitem{alikhani-2010}
  S.~Alikhani, Y.H.~Peng,
  Dominating sets and domination polynomials of certain graphs, II,
  Opuscula Math.\ 30 (2010) 37--51.

\bibitem{alikhani-2014}
  S.~Alikhani, Y.H.~Peng,
  Introduction to domination polynomial of a graph,
  Ars Combin.\ 114 (2014) 257–-266.

\bibitem{anthony-2015}
  B.M.~Anthony, M.E.~Picollelli,
  Complete $r$-partite graphs determined by their domination polynomial,
  Graphs Combin.\ 31 (2015) 1993--2002.

\bibitem{blank-1973}
  M.M.~Blank,
  An estimate of the external stability number of a graph without suspended vertices (in Russian),
  Prikl.\ Mat.\ i Programmirovanie 10 (1973) 3--11.

\bibitem{bresar-2010}
  B.~Bre{\v{s}}ar, S.~Klav{\v{z}}ar, D.~F.~Rall,
  Domination game and an imagination strategy,
  SIAM J.\ Discrete Math. 24 (2010) 979--991.

\bibitem{brouwer-2009}
  A.E.~Brouwer, P.~Csorba, A.~Schrijver,
  The number of dominating sets of a finite graph is odd,
  preprint (June 2, 2009), \path{www.win.tue.nl/~aeb/preprints/domin4a.pdf}.

\bibitem{bujtas-2015}
  Cs.~Bujt\'as,
  Domination game on forests,
  Discrete Math.\ 338 (2015) 2220--2228.

\bibitem{bujtas-2015+}
  Cs.~Bujt\'as, S.~Klav{\v{z}}ar,
  Improved upper bounds on the domination number of graphs with minimum degree at least five,
  Graphs Combin. 32 (2016) 511--519.

\bibitem{fricke-2011}
  G.~Fricke, S.M.~Hedetniemi, S.T.~Hedetniemi, K.R.~Hutson,
  $\gamma$-graphs of graphs,
  Discuss.\ Math.\ Graph Theory 31 (2011) 517--531.

\bibitem{haas-2014}
  R.~Haas, K.~Seyffarth,
  The $k$-dominating graph,
  Graphs Combin.\ 30 (2014) 609–-617.

\bibitem{haddadan-2015}
  A.~Haddadan, T.~Ito, A.E.~Mouawad, N.~Nishimura, H.~Ono, A.~Suzuki, Y.Tebbal,
  The complexity of dominating set reconfiguration,
  Lecture Notes in Comput.\ Sci.\  9214 (2015) 398--409. 


\bibitem{haynes-1998}
  T.W.~Haynes, S.T.~Hedetniemi, P.J.~Slater,
  Fundamentals of Domination in Graphs,
  Marcel Dekker, New York, 1998.
  
\bibitem{mccuaig-1989}
  W.~McCuaig, B.~Shepherd,
  Domination in graphs with minimum degree two,
  J.\ Graph Theory 13 (1989) 749--762.

\bibitem{ore-1962}
  O.~Ore,
  Theory of Graphs,
  American Mathematical Society, Providence, R.I., 1962.

\bibitem{reed-1996}
  B.~Reed,
  Paths, stars and the number three,
  Combin.\ Probab.\ Comput.\ 5 (1996) 277--295.

\bibitem{sohn-2009}
  M.Y.~Sohn, Y.~Xudong,
  Domination in graphs of minimum degree four,
  J.\ Korean Math.\ Soc.\ 46 (2009) 759--773.

\bibitem{suzuki-2014}
  A.~Suzuki, A.E.~Mouawad, N.~Nishimura, 
  Reconfiguration of dominating sets,
  J.\ Comb. Optim., to appear. DOI 10.1007/s10878-015-9947-x.



\end{thebibliography}
\end{document}